\documentclass[12pt]{article}
\usepackage{amsmath,amsthm,amsfonts,amssymb,latexsym}
\date{}

\newtheorem{theorem}{Theorem}

\newtheorem{lemma}[theorem]{Lemma}

\newtheorem{corollary}[theorem]{Corollary}
\newtheorem{observation}[theorem]{Observation}
\newtheorem{proposition}[theorem]{Proposition}

\newtheorem{example}[theorem]{Example}

\newcommand{\IN}{{\mathbb N}}
\newcommand{\IR}{{\mathbb{R}}}
\newcommand{\w}{\omega}
\newcommand{\Te}{{\Theta}}
\newcommand{\Op}{\mathcal O}
\newcommand{\U}{\mathcal U}
\newcommand{\ls}{L_s}

\title{$o$-Boundedness of
free  objects over a Tychonoff space}
\author{Lyubomyr Zdomskyy}
\begin{document}
\maketitle

\begin{abstract}
In this paper we characterize  various sorts of boundedness of the free  (abelian)
topological group $F(X)$ ($A(X)$) as well as the free locally-convex linear topological space
$L(X)$ in terms of properties of a Tychonoff space $X$. 
These properties appear
to be close to so-called selection principles, which permits us to
show, that (it is consistent with ZFC that) the property of Hurewicz (Menger) is $l$-invariant.  
This gives a method of construction of
$OF$-undetermined topological groups with strong combinatorial properties.
\end{abstract}

\section*{Introduction}
\setcounter{section}{1}
\subsection*{Main objects and related notions}

The starting impulse for writing this paper came from \cite{HRT}, where the
problem of characterization of Tychonoff spaces $X$ whose
free (abelian) topological group $F(X)$ ($A(X)$) is [strictly] $o$-bounded was posed.  
In fact, this problem consists of four subproblems. Three of them (except for the characterization
of $o$-boundedness of $F(X)$) are solved here.
Throughout the paper group operations on abelian groups are denoted by $+$
and ``topological space'' means `` Tychonoff space''.

\footnotetext{ \emph{Keywords and phrases.} Selection principle,
(strict) $o$-boundedness, free topological group, $l$-equivalence.
 
\emph{2000 MSC.} Primary: 54D20. Secondary: 22A05, 20E05, 91A44.}

Thus the main objects considered in this paper are free (abelian) topological 
groups over a  space $X$, i.e. a (abelian) topological group $G$ 
that contains $X$ as a set of generators and satisfies the following condition:
each continuous function $\varphi:X\to H$ of $X$ to an arbitrary (abelian) topological group
$H$ admits a unique extension  to a continuous homomorphism $\tilde{\varphi}:G\to H$,
see \cite{GZ} or \cite{Tk2} for basic properties of free topological groups. 
As usually we denote by $C_p(X)$ the space of continuous
real-valued functions on $X$, endowed with topology inherited from the  
Tychonoff product $\IR^X$. It is well-known  \cite[Ch.~0]{Ar89} that the correspondence
$x\mapsto\psi_x$, where $\psi_x(f)=f(x)$ for all $f\in C_p(X)$, is a closed embedding
of $X$ into $C_pC_p(X)$ such that the image of $X$ is linearly independent.
In what follows we denote by $L_p(X)$ the linear hull of $X$ in $C_pC_p(X)$
with the subspace topology. 
The space $L_p(X)$ is  the   free object over $X$ in the category of
linear topological spaces with the weak topology, see \cite{Usp}.
The free object over $X$ in the category of all (locally-convex) linear topological spaces
will be denoted by $\ls(X)$ (resp. $L(X)$).  (The topology on $\ls(X)$
is the strongest linear topology inducing the original topology on the space $X\subset\ls(X)$.
This justifies the choice of the notation for $\ls(X)$.)

 The spaces $X$ and $Y$ are called
\begin{itemize}
\item  \emph{$M$-equivalent},
if the topological groups $F(X)$ and $F(Y)$ are topologically isomorphic;
\item  \emph{$A$-equivalent},
if the topological groups $A(X)$ and $A(Y)$ are topologically isomorphic;
\item \emph{$l$-equivalent}, if $C_p(X)$ and $C_p(Y)$ are isomorphic as linear topological spaces.  
\end{itemize}  
It was shown in \cite{Ar89} that a  space
$X$ is $l$-equivalent to a space $Y$ if and only if $L_p(X)$ is linearly
homeomorphic to $L_p(Y)$. We shall use this as an alternative 
definition of the  $l$-equivalence relation.

  We say that a topological property $\mathsf P$ is $\varphi$-invariant, where $\varphi$
runs over $M$, $A$ and $l$, if a space $X$ has this property whenever so does
any space $Y$ $\varphi$-equivalent  to $X$. It is known \cite{Ar80}
that $M$-equivalence implies $A$-equivalence, and $A$-equivalence implies $l$-equivalence,
and consequently each $l$-invariant property is  $A$-invariant, and each $A$-invariant
property  is $M$-invariant. For various examples of  $\varphi$-invariant
properties see, e.g. \cite[Ch.~2]{Ar89}, \cite{Tk2}, and \cite{Tkachuk}.  In this paper we prove
the $l$-equivalence of selection principles defined below. It is worth to mention
here that these principles are not multiplicative by \cite[Th.~2.12]{JMSS}, which makes
it impossible to use  that the free (nonabelian)  group
over a space $X$ can be represented  as  countable
union of continuous images of  finite powers of $X$.
They are  also not hereditary \cite{BT}, and therefore corresponding proofs can not be
reduced to classical result of V.~Pestov \cite{Pe} asserting that $X$
is a countable union of subspaces each of which is homeomorphic to a subspace
of $Y$ provided $X$ and $Y$ are $M$-equivalent.

\subsection*{Selection principles on topological spaces and groups}

The notion of a  $o$-bounded topological group
was introduced by O.~Okunev and M.~Tka\-chen\-ko with the purpose of characterizing subgroups of
$\sigma$-compact groups, see \cite{He00}, \cite{HRT},  and \cite{Tk98}
for the discussion of these properties. Recall,
that a topological group $G$ is  \emph{$o$-bounded}, if for every
sequence $(U_n)_{n\in\w}$ of nonempty open subsets of $G$, there exists a sequence
$(F_n)_{n\in\w}$ of finite subsets of $G$ such that
$G=\bigcup_{n\in\w}F_n\cdot U_n$. It is clear, that every $\sigma$-compact
group is $o$-bounded.  
Properties of topological spaces $X$ appearing  as duals of the
$o$-boundedness  of free groups and $L_p(X)$ are closely related to so-called selection principles.
The duality between properties of $X$ and $F(X)$, $A(X)$, and $C_p(X)$ is represented
by many results, see \cite{Tk2} and \cite{Ar89}.
The oldests of selection principles, namely the covering properties of Menger 
and Hurewicz\footnote{The
properties $\bigcup_{\mathrm{fin}}(\mathcal O,\mathcal O)$ and
$\bigcup_{\mathrm{fin}}(\mathcal O,\Gamma)$
 in terms of M.~Scheepers \cite{Sch}.},
were introduced at the beginning of 20-th century,
see \cite{Sch} or \cite{Ts-sur-p} for their history and basic properties.
In nearly seventy years after their appearence
M.~Scheepers systematized existing and introduced new properties of this kind.
Among them the property $\bigcup_{\mathrm{fin}}(\mathcal O,\Omega)$ is of extreme
importance for us,  and we shall refer to it as the Scheepers property.
To define the above three selection principles used in this paper, we have to recall
from \cite{GN}  definitions of some classes of covers:
 family $\{U_n:n\in\w\}$ of subsets of $X$ is said
to be
\begin{itemize}
\item an \emph{$\w$-cover} of $X$, if for every finite subset $F$ of $X$ there exists
$n\in\w$ such that $F\subset U_n$; 
\item a \emph{$\gamma$-cover} of $X$,  if for every $x\in X$ the set $\{n\in\w:x\not\in U_n\}$
is finite.
\end{itemize}
Let $B$ be a subset of a set $X$ and $u$ be  a cover of $X$. We say that $B$ is
\emph{$u$-bounded}, if $B\subset\cup c$ for some finite subfamily $c$ of $u$.  
A topological space $X$ is said to have \emph{the Menger} (resp. \emph{Scheepers, Hurewicz})
property, if  for every sequence $(u_n)_{n\in\w}$  of open covers of $X$ 
there exists a ($\w$-, $\gamma$-) cover $\{B_n:n\in\w\}$ of $X$ such that each $B_n$ is $u_n$-bounded.
The Menger and Scheepers properties differ under the Continuum Hypothesis
by \cite[Theorem~2.8]{JMSS} and coincide under $(\mathfrak u<\mathfrak g)$ according to
Corollary~2 of \cite{Zd1}.
Note that the $o$-boundedness (of all finite powers) of a topological group is nothing else but the Menger
(Scheepers) property applied to the family of uniform covers with 
respect to the left uniformity on it, see Lemma~\ref{l2} and Proposition~\ref{basics}. 
From now on we denote by $\nu X$ and $\mu X$ the Hewitt and Dieudonne completions
of a space $X$, see \cite{En} for their definitions and basic properties.
We are in a position now
to present the characterization of the $o$-boundedness of free abelian topological group.
 
\begin{theorem} \label{main1} 
For a space $X$ the following conditions are equivalent:
\begin{itemize}
\item[$(1)$]  $A(X)$ is $o$-bounded;
\item[$(2)$] $A(X)^n$ is $o$-bounded for all $n\in\IN$;
\item[$(3)$] $L_p(X)$ is $o$-bounded;
\item[$(4)$] $L(X)$ is $o$-bounded;
\item[$(5)$] $\ls(X)^n$ is $o$-bounded for all $n\in\IN$;
\item[$(6)$] every continuous metrizable image of $X$ has the Scheepers property;
\item[$(7)$] $A(\nu X)$ is $o$-bounded;
\item[$(8)$] $A(\mu X)$ is $o$-bounded.
\end{itemize} 
\end{theorem}

\subsection*{Selection games and multicovered spaces}

$o$-Boundedness as well as the Menger property have natural game counterparts.
In case of a $\sigma$-compact group
$G$ a sequence $(F_n)_{n\in\w}$  witnessing the $o$-boundedness of $G$
may be constructed
by the second player in the process of an infinite game, called $OF$.
This game is played by two players, say I and II. Player I selects an open subset
$U_0$ of $G$, and player two responds choosing  some finite subset $F_0$ of $G$. 
In the second turn, player $I$ selects some open subset $U_1$ of $G$,
and II responds choosing a finite subset $F_2$ of $G$, and so on. At the end of this game
we obtain  the  sequences $(U_n)_{n\in\w}$ and $(F_n)_{n\in\w}$.
Player II is declared the winner, if $\bigcup_{n\in\w}F_n\cdot U_n=G$. Otherwise,
player I wins. A group $G$ is \emph{strictly $o$-bounded}, if the second player
(= player II) has a winning strategy in the game $OF$ on $G$. 
If none of the players has a winning strategy, then  $G $ is called
\emph{$OF$-undetermined}.
It is clear that each $\sigma$-compact group is strictly $o$-bounded and
thus $o$-bounded. Examples distinguishing the $\sigma$-compacness, strict
$o$-boundedness and $o$-boundedness may be found
in  \cite{BNS}, \cite{He00}, \cite{Tk98},     and \cite{Ts}.

As we shall see later,  the strict $o$-boundedness
of the free objects over a space $X$  has no characterization in terms
of continuous metrizable images of $X$ in spirit of Theorem~\ref{main1}.
This  constrained us to use the language of multicovered spaces, which seems 
to be the most appropriate one for description of the corresponding
 property of $X$. 
By a \emph{multicovered space}\footnote{
The notion of a multicovered space and some other notions related to multicovered spaces,
as well as Corollaries~\ref{cb1} and \ref{Banakh} are due to T.Banakh. Multicovered spaces, which seem to be
the most general objects where properties looking similar to (strict) $o$-boundedness
can be considered, are discussed in \cite{BZ}.} we understand a pair
$(X,\lambda)$, where $X$ is a set and  $\lambda$  is a \emph{multicover} of
$X$, i.e.  a family of covers of $X$. 
There are many natural examples of multicovered spaces:
\begin{itemize}
\item Each topological space $ X $ can be considered as  a multicovered space
$ (X,\Op)$, where $\Op$ denotes the family of all open covers of
$X$;
\item Every metric space $(X,\rho)$ admits a natural  multicover $\lambda_\rho$
consisting of covers by $\varepsilon$-balls: $\lambda_\rho=\{\{B_\rho(x,\varepsilon):x\in X\}:\varepsilon>0\}$,
where $B_\rho(x,\varepsilon) = \{y\in X:\rho(y,x)<\varepsilon\}$;
\item Every uniform space $(X,\mathcal{U})$ has a multicover $\lambda_\mathcal{U}$ consisting of uniform covers,
i.e. $\lambda_\mathcal{U}=\{\{U(x):x\in X\}:U\in\mathcal{U}\}$, where $U(x)=\{y\in X:(x,y)\in U\}$;
\item In particular, each topological group $G$ admits four natural multicovers
$\lambda_L(G)$, $\lambda_R(G)$, $\lambda_{L\vee R}(G)$ and $\lambda_{L\wedge R}(G)$ corresponding to its
left, right, two-side and R\"{o}lcke uniformities, see \cite{RD} for
more information on these uniformities;   
\item In case of an abelian topological group $G$ all of the above
uniformities coincide, and we denote them by $\U(G)$. The family
$\{\{(x,y):x-y\in U\}:0\in U \in\mathcal O(G)\}$ is a base of $\U(G)$.
Therefore corresponding multicovers coincide as well, and we denote them by
$\lambda(G)$. 
\end{itemize} 

By analogy with the game $OF$ on a topological group $G$ we can introduce the game
$CB$ (abbreviated from Cover-Bounded) on a multicovered space $(X,\lambda)$ as follows:
two players, I and II, step by step choose a cover $u_n\in\lambda$ and an
$u_n$-bounded subset $B_n$ of $X$, respectively. The player II is declared the winner,
if $X=\bigcup_{n\in\w}B_n$. Otherwise the player I wins.
A multicovered space $(X,\lambda)$ is said to be \emph{winning},
if the second player has a winning strategy in the game $CB$ on $(X,\lambda)$.
It is clear that that the game $OF$ on a topological group $G$ is
equivalent to the game $CB$ on the multicovered space $(G,\lambda_L)$
in the sense that one of the players has a winning strategy in one of these games
if and only if he has a winning strategy in the other one.
It also should be mentioned here that the game $CB$ on a multicovered space
$(X,\mathcal O(X))$ is nothing else but the game $H(X)$
  introduced by R.~Telgarsky in \cite{Topsoe}, see also  \cite{Te2}
and references there in. 

Let $X$ be a Tychonoff space. Recall from \cite{En} that the  uniformity
$\U$ on $X$ is called \emph{universal}, if it generates the topology
of $X$ and contains all uniformities on $X$ with this property.
Throughout this paper the universal uniformity of a topological space
$X$ will be denoted by $\U(X)$. The reader is refered to the next section for the definition
of the product of multicovered spaces.
We are in a position now to present the main result of this paper.

\begin{theorem} \label{main-win} 
For a space $X$ the following conditions are equivalent:
\begin{itemize}
\item[$(1)$] $(F(X),\lambda_{L\wedge R})$ is winning;
\item[$(2)$] $F(X)$ is strictly $o$-bounded;
\item[$(3)$] $(F(X)^n,\lambda^n_{L\vee R})$ is winning
for all $n\in\w$;
\item[$(4)$]  $A(X)$ is strictly $o$-bounded;
\item[$(5)$]  $L_p(X)$  is strictly $o$-bounded;
\item[$(6)$] $L(X)$ is strictly $o$-bounded;
\item[$(7)$]  $\ls(X)^n$  is strictly $o$-bounded for all
$n\in\IN$;
 \item[$(8)$] $(X,\lambda_{\U(X)})$ is winning.
\end{itemize} 
\end{theorem}

The equivalent properties of in the above theorem are not preserved by finite powers.
To describe a corresponding space we have to introduce some  notions related
to multicovered spaces. 
A multicovered space $(X,\lambda)$ is called
\begin{itemize}
\item \emph{totally-bounded}, if $X$ is  $u$-bounded for every $u\in\lambda$; 
\item \emph{$\w$-bounded}, if each cover $u\in\lambda$ has a countable subcover.
\end{itemize}

These notions generalize the $\w$-boundedness
of uniform spaces introduced by I.~Guran in \cite{Gu} and the total boundedness 
in sense that a uniform space $(X,\U)$ has one of the above properties if and only
if so does the multicovered space $(X,\lambda_{\U})$. For example, $(X,\Op(X))$
is totally-bounded ($\w$-bounded) if and only if $X$ is compact (Lindel\"of).
Let $X$ be a countably-compact spaces $X$
such that there exists a continuous pseudometric $\rho$ on $X^2$
such that the space $X^2$ is not Lindel\"of, see Example~\ref{count-com}. Then
the uniform space $(X,\U(X))$ as well as  the multicovered space $(X,\lambda_{\U(X)})$
are totally-bounded, and consequently $(X,\lambda_{\U(X)})$ is winning.
But the uniform space $(X^2,\U(X^2))$ obviously fails to be $\w$-bounded,
consequently so does the multicovered space $(X^2,\lambda_{\U(X^2)})$,
and hence  $X^2$ does not satisfy the conditions  of Theorem~\ref{main-win}.

The properties of Menger, Scheepers, and Hurewicz can be also naturally carried out in the realm of multicovered spaces:
a multicovered space $(X,\lambda)$ has the Menger (resp. Scheepers, Hurewicz)
property if for every sequence $(u_n)_{n\in\w}\in\lambda^\w$ there exists a
sequence $(B_n)_{n\in\w}$ of subsets of $X$ such that each $B_n$
is $u_n$-bounded and $\{B_n:n\in\w\}$ is a cover (resp. $\w$-cover, $\gamma$-cover)
of $X$. It is a simple exercise to show that each Menger multicovered space
is $\w$-bounded.   A crucial observation here is that a topological group $G$ is $o$-bounded if and only
if the multicovered space $(G,\lambda_R(G))$ is Menger. The $o$-boundedness of
free objects may be also described in terms of properties 
of the multicovered space $(X,\lambda_{\U(X)})$ as well,
which extends Theorem~\ref{main1}.

\begin{theorem} \label{main-mult1} 
Let $X$ be a Tychonoff space. Then $A(X)$ is $o$-bounded
if and only if the multicovered space $(X,\lambda_{\U(X)})$ is Scheepers.
\end{theorem}

The Hurewicz property is selfdual.

\begin{theorem}  \label{main-hur}  
For  a  space $X$ the following conditions are equivalent:
\begin{itemize}
\item[$(1)$]  $(A(X),\lambda(A(X)))$ is Hurewicz;
\item[$(2)$]  $(A(X)^n,\lambda(A(X))^n)$ is Hurewicz for all $n\in\w$;
\item[$(3)$] $(L_p(X),\lambda(L_p(X)))$ is Hurewicz;
\item[$(4)$] $(L(X),\lambda(L(X)))$ is Hurewicz;
\item[$(5)$] $\ls(X)^n$ is Hurewicz for all $n\in\IN$;
\item[$(6)$]  every continuous metrizable image of $X$ is Hurewicz;
\item[$(7)$] $(A(X),\lambda(A(\nu X)))$ is Hurewicz;
\item[$(8)$] $(A(X),\lambda(A(\mu X)))$ is Hurewicz.
\item[$(9)$] $(X,\lambda_{\U(X)})$ is Hurewicz.
\end{itemize} 
\end{theorem}

For a Lindel\"of topological space $X$ the  multicovers $\lambda_{\U(X)}$ and 
$\Op(X)$ are equivalent in the sense defined in the next section, see
Corollary~\ref{cle213}. In combination with Proposition~\ref{basics}
this gives us  the subsequent
\begin{proposition} \label{Lind-top} 
The multicovered space $(X,\Op(X))$ is winning (resp. Menger, Scheepers, Hurewicz)
if and only if so is $(X,\lambda_{\U(X)})$ and $X$ is Lindel\"of.
\end{proposition}
Note, that the multicovered space $(X,\Op(X))$ has the Menger (resp. Scheepers, Hurewicz) property
if and only if so does the topological space $X$.
Concerning the winning property of the multicovered space $(X,\Op(X))$, there are
 many equivalent  statements to it. At the beginning of 80-th
R.~Telgarsky introduced the game $H(X)$ 
(implicitly existing in earlier works of W.~Hurewicz)
on a topological space $X$, which coincides with the game $OF$
on the multicovered space $(X,\Op(X))$, and proved that the second player has a winning strategy
in this game if and only if $X$ is $C$-like, which means that the first player has a winning strategy in the compact-open game
on $X$, see \cite{Te2} and references there in. 
In the current terminology the game $H(X)$ is called the Menger game on $X$,
see \cite{Sch}.
This yields the subsequent reformulation of Proposition~\ref{Lind-top}.
\begin{corollary} \label{c11} 
A Tychonoff space $X$ is $C$-like
(resp. Menger, Scheepers, Hurewicz)
if and only if it is Lindel\"of and the multicovered space
$(X,\lambda_{\U(X)})$ is winning (resp. Menger, Scheepers, Hurewicz).
\end{corollary}
Since the Lindel\"of property is $l$-invariant \cite{Velichko}
(see also \cite{Bouziad}, where it is shown that the Lindel\"of number is $l$-invariant),
and the properties of the multicovered space $(X,\lambda_{\U(X)})$ considered in the
above corollary  have counterparts
among the properties of $L_p(X)$ obviously preserved by linear homeomorphisms,
we get the following

\begin{corollary} \label{invariant} 
The properties of Scheepers, Hurewicz and being $C$-like are
$l$-invariant (and hence $A$- and $M$-invariant).
Consequently the Menger property is $l$-invariant under $\mathfrak u<\mathfrak g$
being equivalent to the Scheepers one.
\end{corollary}

In light of this it is worth to mention Question~II.2.8  of \cite{Ar89}
whether the Menger property is $t$-invariant, which seems to
be still unsolved.

As we could see in Theorems~\ref{main1},
\ref{main-mult1}, and \ref{main-hur}, the Scheepers and Hurewicz properties
of the multicovered space $(X,\lambda_{\U(X)})$ admit a characterization in terms
of metrizable images of $X$. In light of this 
 one may try to prove some similar characterization of the winning property
of  $(X,\lambda_{\U(X)})$ of the following kind: $(X,\lambda_{\U(X)})$ is winning
if and only if every metrizable image of $X$ has some ``strong'' property $\mathsf P$.
But even the property $\mathsf P$ of being countable, which seems to be the strongest
among those one could consider, does not work.
Let us recall that a topological space $X$ is a \emph{$P$-space}, if  every $G_\delta$-subset
of $X$ is open.
R.~Telgarsky in \cite{Te4} observed that 
  a Lindel\"of $P$-space $Y$ constructed by R.~Pol in \cite{Pol} 
 fails to be $C$-like (and hence $(X,\lambda_{\U(X)})$ is not winning by Corollary~\ref{c11}).
It sufficies to note that the size of arbitrary metrizable  image of a $P$-space $X$
does   not exceed the Lindel\"of number of $X$. 

\begin{corollary} \label{Tel-Pol} 
The Pol's space $Y$ has the following
properties: 
\begin{itemize}
\item[$(1)$] The groups $F(Y)$, $A(Y)$, $L_p(Y)$, $L(X)$, and $\ls(X)$
are $OF$-undetermined;  
\item[$(2)$] all metrizable images of $Y^n$, $F(Y)$, $A(Y)$ are countable,
where $n\in\IN$.
\end{itemize}
\end{corollary}
\begin{proof}
The second item obviously
follows from the facts that every finite power of a Lindel\"of $P$-space
is again Lindel\"of $P$-space, and each metrizable image of a Lindel\"of
$P$-space is countable. 
 Concerning the first one,
it simply follows  from Corollary~\ref{afterend} and an observation that
each Lindel\"of $P$-space is Hurewicz.
\end{proof}
Corollary~\ref{Tel-Pol} is closely related to the result of A.~Krawczyk
and H.~Michalewski \cite{KM}
who used the space $Y$ to construct an $OF$-undetermined  $P$-group $G$.
Similar ideas are also used in Theorem~3.1 of \cite{HRT}.

The problem of construction of $OF$-undetermined groups was posed in \cite{Tk98}
and solved in \cite{KM} and \cite{Ba02} (and, probably, somewhere else)
independently. Theorems~\ref{main-win} and \ref{main-hur} 
supply us with a method of constructing 
$OF$-undetermined groups:  it sufficies to take a topological space
which does not satisfy condition $(8)$ of Theorem~\ref{main-win}
but still has some strong property guaranteeing the first player
having no winning strategy on free objects considered in this paper.
Such the  properties are given by the subsequent proposition,
which easily follows from \cite[Theorem 27]{JMSS}.
\begin{proposition} \label{neznaju} 
Let $G$ be a topological group such that the underlying topological space is Menger.
Then the first player has no winning strategy in the game $OF$
on $G$.
\end{proposition}

Finally, we present (nonmetrizable) examples of  $OF$-undetermined
groups.

\begin{corollary} \label{afterend} 
Let $X$ be a non-$\sigma$-compact  metrizable space such that  all finite
powers of $X$ are Menger (Hurewicz). Then all finite powers of $F(X)$, $A(X)$, $L_p(X)$, $L(X)$, and $\ls(X)$
are $OF$-undetermined being Menger (Hurewicz) groups which fail to be
strictly $o$-bounded.
\end{corollary}
\begin{proof}
Let $X$ be a non-$\sigma$-compact metrizable  space whose all finite
powers of $X$ are Hurewicz (Menger),
 and $G$ be one of the  groups $F(X)$, $A(X)$, $L_p(X)$, $L(X)$, $\ls(X)$.
Since $G$ is the countable union of continuous images of finite powers of $Y=X\times\IR$
(see the proof of Theorem~\ref{main1}, where it is  shown that $\ls(X)$
is a contiuous image of $A(X\times \IR)$), so is $G^n$ for all $n\in\IN$.
 Applying   Lemma~\ref{l5} and Corollary~\ref{cle213},  we conclude
that  each finite power of 
 of $Y$ is Hurewicz
(Menger).
 Since the Hurewicz (Menger) property is preserved by
continuous  images and countable unions by Lemma~\ref{l6}
(this was also pointed out in \cite{JMSS}),
$G^n$ is Hurewicz (Menger) for all $n\in\w$. Applying Proposition~\ref{neznaju},
we conclude that for every $n\in\w$ the first player has no winning strategy in the game $OF$
on $G^n$.

As it was shown by R. Telgarsky, every winning (= $C$-like) metrizable topological space
is $\sigma$-compact, see \cite{direct}. Therefore $X$ fails to be winning,
and thus $G$ is not strictly $o$-bounded by Theorem~\ref{main-win}. Consequently
$G^n$ is not strictly $o$-bounded for all $n\in\IN$.
From the above it follows that $G^n$ is $OF$-undetermined for all $n\in\IN$.
\end{proof} 
\begin{observation}
Every topological group 
with the Hurewicz property is strictly $o$-bounded provided it is metrizable.
\end{observation}
\begin{proof}
Let $G$ be a  topological group whose underlying  topological space $G$ is Hurewicz 
and $\{U_n:n\in\w\}$
be a countable local base at the identity of $G$. Set $u_n=\{gU_n:g\in G\}$.
The Hurewicz property of $G$ yields a sequence $(F_n)_{n\in\w}$ 
of finite subsets of $G$ such that 
$G=\bigcup_{n\in\w}\bigcap_{k\geq n}F_n U_n$,
consequently $G$ is strictly $o$-bounded as a countable union of its totally-bounded
subspaces.  
\end{proof}

Spaces $X$ with such propertis as in  Corollary~\ref{afterend} were constructed in
\cite{BT}, \cite{CP}, and \cite{TZ}.

\section*{Proofs}
\setcounter{section}{2}
In our proofs of results announced in Introduction
we shall exploit a number of auxiliary statements about multicovered
spaces. As a matter of the fact,
all of these results are (consequences of more general ones) proven in \cite{BZ}.
But in sake of completeness we present their proofs.
 Their formulations involve some additional notions and notations.
For  uniform spaces $(X_1,\U_1)$ and $(X_2,\U_2)$ we shall 
identify
the uniformity on their product $X_1\times X_2$ generated by $\U_1$ and $\U_2$
with the product $\U_1\times\U_2$. 
Let $u$ and $\lambda$ be a cover and a multicover of a set $X$ respectively, and $Z\subset X$. Then
$u|Z$ denotes the family $\{U\cap Z:U\in u\}$ and  $\lambda|Z=\{u|Z:u\in\lambda\}$.
Every subset $Z$ of $X$ with the induced multicover $\lambda|Z$
is called \emph{a subspace} of the multicovered space $(X,\lambda)$.
By the product  of multicovers $\lambda$ and $\nu$ of sets $X$ and $Y$
we understand the multicover $\eta=\{u\cdot v:u\in\lambda,v\in\nu\}$ of $X\times Y$, 
where $u\cdot v=\{U\times V:U\in u, V\in v\}$. 
Again, we identify $\eta$ with the product $\lambda\times\nu$.

Next, we shall also use the  preorder $\prec$ on the family of all covers of a set
$X$, where $u\prec v$ means that each $v$-bounded subset is $u$-bounded. In other words,
$u\prec v$ if and only if for every finite subset $c$ of $v$ there exists a finite subset
$d$ of $u$ with $\cup d\supset\cup c$.
Note, that all multicovers $\lambda$ considered in this paper are \emph{centered},
which means that each finite subset $c$ of $\lambda$ has an upper bound in $\lambda$
with respect to $\prec$.  Let us also observe that $u\prec v$ provided $v$
is a refinement of $u$ in the sense that each $V\in v$ lies in some $U\in u$.
 The preorder $\prec$ on the family of all
covers of $X$ generates the following preorder on the family of all multicovers of $X$,
which is also denoted by $\prec$: 
$\lambda\prec \nu$ if and only if for every $u\in \lambda$ there exists $v\in\nu$
such that $u\prec v$. We say that multicovers $\lambda$ and $\nu$ of a set $X$
are \emph{equivalent} (and  write $\lambda\cong\nu$) if $\lambda\prec\nu$ and $\nu\prec\lambda$.
Given any multicovers $\lambda$ and $\nu$ of $X$ and $Y$ respectively, we call a function $f:X\rightarrow Y$
\begin{itemize}
\item \emph{uniformly bounded},
if for every $v\in\nu$ there exists  $u\in\lambda$ such that for every $u$-bounded 
subset $A$ of $X$ its image $f(A)$ is $v$-bounded;
\item  \emph{perfect},
if for every $u\in\lambda$  there exists  $v\in\nu$  such that for every $v$-bounded 
subset $B$ of $Y$ the preimage $f^{-1}(B)$ is $u$-bounded.
\end{itemize}

In the subsequent simple statement we collect some straightforward
properties of the notions introduced before.
\begin{proposition} \label{basics} 
\begin{itemize}
\item[$(1)$] Let $(X_1,\U_1)$ and $(X_2,\U_2)$ be uniform spaces. Then the multicovers
$\lambda_{\U_1\times\U_2}$ and $\lambda_{\U_1}\times\lambda_{\U_2}$ of $X_1\times X_2$
are equivalent.
\item[$(2)$] Let $G$ and $H$ be  topological groups and $T\in\{L,R,L\vee R, L\wedge R\}$. Then $\lambda_{T}(G\times H)$ 
is equivalent to $\lambda_T(G)\times\lambda_T(H)$.
\item[$(3)$] Let $G$ be a  topological group $G$. Then $\lambda_L(G)\cong\{\{gU:g\in
G\}:e\in U\in\Op(G)\}$, \\
 $\lambda_R(G)\cong\{\{Ug:g\in G\}:e\in U\in\Op(G)\}$, \\
 $\lambda_{L\vee R}(G)\cong\{\{gU\cap Ug:g\in G\}:e\in U\in\Op(G)\}$, \\
$\lambda_{L\wedge R}(G)\cong\{\{UgU:g\in G\}:e\in U\in\Op(G)\}$.
\item[$(3')$] For an abelian group $G$ the multicover $\lambda(G)$
is equivalent to $\{\{g+U:g\in G\}:e\in U\in\Op(G)\}$.
\item[$(4)$] Let $\lambda$ and $\nu$ be multicovers of a set $X$ and $\lambda\prec\nu$.
Then $(X,\lambda)$ is Menger (resp. Scheepers, Hurewicz, winning) provided so is $(X,\nu)$.
\item[$(5)$] 	If  multicovers $\lambda$ and $\nu$ of $X$ are equivalent,
then $(X,\lambda)$ is Menger (resp. Scheepers, Hurewicz, winning) if and only if so is $(X,\nu)$.
\item[$(6)$] $\lambda\prec\nu$ if and only if the identity map $\mathrm{id}_X$ is perfect
with respect to $\lambda$ and $\nu$.
\item[$(7)$] If $f:X\rightarrow Y $ is perfect with respect
to multicovers $\lambda$ and $\nu$ of $X$ and $Y$ respectively, then 
 $(X,\lambda)$ is Menger (resp. Scheepers, Hurewicz, winning) provided so is $(Y,\nu)$.
\item[$(8)$] If $f:X\rightarrow Y$ is uniformly bounded with respect to
multicovers $\lambda$ and $\nu$, then $(Y,\nu)$ is winning (resp. Hurewicz,
Scheepers, Menger) provided so is $(X,\lambda)$. 
\end{itemize}
\end{proposition}
\begin{proof}
Because of simplicity of all of the items, we shall only proof the ``winning'' part
of the seventh  one. For this purpose we have to consider more formally
the notion of a winning strategy in the game $CB$ on a multicovered space $(X,\lambda)$. 
By a strategy of a  second player we understand a map $\Theta:\lambda^{<\omega}\rightarrow\mathcal{P}(X)$
assigning to each finite sequence of covers $(u_0,\ldots,u_n)\in\lambda^{<\omega}$
a $u_n$-bounded subset $\Theta(u_0,\ldots,u_n)$ of $X$, where $\lambda^{<\w}=\bigcup_{n\in\w}\lambda^n$. 
Such a strategy is \emph{winning},
if the family $\{\Theta(u_0,\ldots,u_n):n\in\omega\}$ is a cover of $X$
 for any  sequence $(u_n)_{n\in\omega}$.

Now, assume that $\Te_Y$ is a winning strategy of the second player in the game $CB$
on a multicovered space $(Y,\nu)$. Construct a map $\phi:\lambda\rightarrow\nu$
such that $f^{-1}(B)$ is $u$-bounded for every $\phi(u)$-bounded subset 
$B$ of $Y$. It sufficies to observe that 
 $\Te_X:(u_0,\ldots,u_n)\mapsto f^{-1}(\Te_Y(\phi(u_0),\ldots,\phi(u_n)))$
is a winning strategy of the second player in the game $CB$ on $(X,\lambda)$.
\end{proof}

 We shall exploit the following important result of V.~Pestov, see \cite{Pea} or \cite[2.8]{Tk2}.
\begin{proposition} \label{fromG} 
Let $X$ be a Tychonoff space. Then  the 
natural uniformity on $A(X)$ generates the universal uniformity $\U(X)$ on $X$.
\end{proposition}

We shall also use the following straightforward consequence of 
Lemma~1.0 of \cite{Le}.


\begin{corollary}\label{le213} 
If $X$ is a Lindel\"{o}f regular space and $u$ is an 
 open cover of $X$, then there exists a pseudometric $d$ on $X$
 such that a subset $Y$ of $X$ is $u$-bounded 
 provided $\mathrm{diam}_d(Y)<\infty$
(Here, as usual, $\mathrm{diam}_d(Y)=\sup\{d(y_1,y_2):y_1,y_2\in Y\}$).
\end{corollary}

\begin{corollary} \label{cle213} 
Let $X$ be a Lindel\"of regular space. Then the multicovers $\mathcal O(X)$
and $\lambda_{\U(X)}$ are equivalent.
\end{corollary}
\begin{proof}
Since every uniform cover has an open uniform refinement,
we conclude that $\lambda_{\U(X)}\prec\Op(X)$.

To prove that $\Op(X)\prec\lambda_{\U(X)}$, fix an arbitrary open cover $u$
of $X$ and find a pseudometric $d$ on $X$ such as in Lemma~\ref{le213}.
Then for the uniform cover $v=\{B_d(x,1):x\in X\}\in\lambda_{\U(X)}$ 
we obviously have $u\prec v$, which finishes our proof. 
\end{proof}

Let $A$ be a subset of the Cartesian product $X\times Y$. From now on we shall use the following
notations:
$A^{-1}=\{(y,x)\in Y\times X:(x,y)\in A\}$, $A(x)=\{y\in Y:(x,y)\in A\}$, where
$x\in X$. 
Recall from \cite{Gu} that a unifom space $(X,\U)$ is \emph{$\w$-bounded},
if each uniform cover contains a countable subcover. In particular, topological
group $G$ is $\w$-bounded, if so is the uniform space $(G,\U)$, where
$\U$ is the left uniformity of $G$.

\begin{lemma}\label{l1} 
Let $X$ be a  space such that $(X,\U(X))$ is $\w$-bounded and $G\supset X$ be an abelian topological group
such that every continuous map $\phi:X\to\IR$ can be extended to a continuous
homomorphism $\tilde{\phi}:G\to\IR$. 
Then the maps
  $\psi_n:X^n\rightarrow G$, 
where
$\psi_n:(x_1,x_2\ldots,x_n)\mapsto x_1+ x_2+ \cdots + x_n$,
 are perfect with respect to   $\lambda_{\U(X)^n}$
and $\lambda(G)$ for all $n\in\IN$.
\end{lemma}
\begin{proof}
Given any $u_0\cdot\:\cdots\:\cdot u_{n-1}\in\lambda_{\U(X)^n}$,
find $U\in\U(X)$ such that the uniform cover $u=\{U(x):x\in X\}$ is an upper bound
of the family $\{u_i:i<n\} $ with respect to $\prec$. Let $\rho$
be a pseudometric on $X$ such that $\{(x,y)\in X^2:\rho(x,y)<1\}\subset U$.
Since $(X,\U(X))$ is $\w$-bounded, the space $(X,\rho)$ is Lindel\"of.
Applying Lemma~\ref{le213} to the regular Lindel\"of space $(X,\rho)$
and the cover $w=\{B_\rho(x,1):x\in X\}$, we can find a continuous pseudometric
$d$ on $X$ such that each $Y\subset X$ is $w$-bounded provided $\mathrm{diam}_d(Y)<\infty$.
Fix arbitrary $x_0\in X$ and define a map $f:X\rightarrow\IR$ letting
$f(x)=d(x,x_0)$. From the above it follows that $f^{-1}(-r,r)$ is
$w$-bounded, and hence $u$-bounded for every $r\in\IR$.
Let $\hat{f}:G\rightarrow\IR$ be a continuous
homomorphism extending $f$ and $O$ be an open neighborhood
of the identity of $G$ such that 
$\hat{f}(O)\subset (-1,1)$.

Let us fix arbitrary finite subset $K$ of $G$.
Our proof will be completed as soon as we shall show that
$B=\psi_n^{-1}(O+K)$ is $w^n$-bounded.
By our choice of $O$ there exists $r>0$ such that 
$\hat{f}(O+K)\subset(-r,r)$. Therefore,
$B\subset\psi_n^{-1}(\hat{f}^{-1}(-r,r))$. Let us note, that  
$\hat{f}\circ\psi_n(x_1,\ldots,x_n)=f(x_1)+\ldots+f(x_n)$, consequently
$0\leq f(x_i)<r$ for every $(x_1,\ldots,x_n)\in B$ and $i\leq n$,
and finally  $B$ is $w^n$-bounded being a subset
of $(f^{-1}(-r,r))^n$.
\end{proof}

Next, we shall deal with preservation of selection principles by
operations of finite products and countable unions.   

\begin{lemma}\label{l2}  
A multicovered space $(X,\lambda)$ is Scheepers if and only if $(X^n,\lambda^n)$
is Menger for all $n\in\w$. Consequently the class of Scheepers multicovered
spaces is closed under taking finite powers of its elements.
\end{lemma}
\begin{proof}
Suppose that $(X^n,\lambda^{n})$ is Menger for every $ n\in\omega $.
To see that $(X,\lambda)$ is Scheepers, fix any sequence
$(u_n)_{n\in\omega}\in\lambda^\w$. For every $n\in\omega$ we can apply
the Menger property of $X^n$ to find a cover $\{B_{n,k}^n:k\ge
n\}$ of $X^n$ by powers of $u_k$-bounded sets $B_{n,k}\subset X$.

For every $k\in\w$ let $B_k=\bigcup_{n\le k}B_{n,k}$. We claim
that $\{B_n:n\in\w\}$ is an $\omega $-cover of $ X $. Indeed, fix
any finite subset $F=\{x_1,\ldots, x_n\}$ of $X$. Since the family
$\{B^n_{n,k}:k\geq n\}$ covers $X^n$, $(x_1,\ldots,x_n)\in
B^n_{n,k}$ for some $k\ge n$. Consequently, $F\subset
B_{n,k}\subset B_k$, which completes the proof of the ``if'' part.

To prove the ``only if'' part, suppose that $(X,\lambda)$ is Scheepers.
 To show that the powers of $X$ are Menger, fix some
$ n\in\omega $ and a sequence of covers $(w_k)_{k\in\omega}$ in
$\lambda^{n}$. For every $ k\in\omega $ we can write $ w_{k} $ in the
form $ w_k=u_{k1}\cdot\ldots \cdot u_{kn}$, where $ u_{ki}\in\lambda
$, $ i\in\{1,\ldots,n\}$. Since $(X,\lambda)$ is centered,
 for every $ k\in\omega $ we can find $ u_k\in\lambda $ such that
 $u_k\succ u_{ki}$ for all $ i\in\{1,\ldots,n\}$.
Using the Scheepers property of $X$, find an $\w$-cover
$\{B_k:k\in\w\}$ by a $u_k$-bounded subsets $B_k\subset X$. We
claim that $X^n=\bigcup_{k\in\w}B_k^n$, which clearly implies the
Mengerness of $(X^n,\lambda^{n})$. Indeed, fix any $
x=(x_1,\ldots,x_n)\in X^n $ and find $ k\in\omega $ such that
$\{x_1,\ldots,x_n\}\subset B_k$. Then $x\in B_k^n$.
\end{proof}

We need the following auxiliary notion:
a family $\{A_n:n\in\w\}$ is called a \emph{proper $\w$-cover} of a set $X$,
if for every finite subset $K$ of $X$ the set $\{n\in\w:K\subset A_n\}$ is infinite.

\begin{lemma} \label{l3} 
Let $(X,\lambda)$ be a Scheepers multicovered space. Then for each sequence
$(u_n)_{n\in\w}\in\lambda^\w$ there exists a proper $\w$-cover $\{B_n:n\in\w\}$
of $X$ such that $B_n$ is $u_n$-bounded for all $n\in\w$.
\end{lemma}
\begin{proof}
Let $(u_n)_{n\in\omega}\in\lambda^{\omega}$ be a sequence of covers of
$X$. Using the Scheepers property of $(X,\lambda)$, for every
$k\in\omega$ we can find a sequence $(A_{k,n})_{n\ge k}$ of
$u_n$-bounded subsets $A_{k,n}\subset X$ such that the 
family $\{A_{k,n}:n\ge k\}$ is an $\w$-cover of $X$. For every
$n\in\w$ consider the $u_n$-bounded subset $B_n=\bigcup_{k\le
n}A_{k,n}$ of $X$ and note that $\{B_n:n\in\w\}$ is a proper
$\w$-cover of $X$, which finishes our proof. 
\end{proof}

\begin{lemma} \label{lh} 
The product $(X\times Y,\lambda_X\cdot\lambda_Y)$ of Hurewicz multicovered spaces $(X,\lambda_X)$
and $(Y,\lambda_Y)$ is Hurewicz. Consequently the class of Hurewicz multicovered
spaces is closed under taking finite products of its elements.
\end{lemma}
\begin{proof}
Let us fix a sequence $(w_n)_{n\in\w}\in (\lambda_X\cdot\lambda_Y)^\w$. For every 
$n\in\w$ find $u_n\in\lambda_X$ and $v_n\in\lambda_Y$ such that $w_n=u_n\cdot v_n$.
By the definition of the Hurewicz property,
there are sequences $(A_n)_{n\in\w}$ and $(B_n)_{n\in\w}$ of subsets of $X$ and $Y$
respectively such that each $A_n$ ($B_n$) is $u_n$- ($v_n$-) bounded,
and the families $\{A_n:n\in\w\}$  and
$\{B_n:n\in\w\}$ are $\gamma$-covers of $Y$. 
 For every $n\in\w$ put $C_n=A_n\times B_n$. It is a simple matter to verify that
the family $\{C_n:n\in\w\}$ is a $\gamma$-cover  of $X\times Y$ and each $C_n$ is $w_n$-bounded,
which finishes our proof.
\end{proof}

\begin{lemma} \label{l6} 
 Let $A_n$, $n\in\omega$, be subspaces of a
multicovered space $(X,\lambda)$. If every subspace $A_n$, $n\in\omega$,
is winning (resp. Menger, Hurewicz), then so is their union
$A=\bigcup_{n\in\omega}A_n$.
\end{lemma}
\begin{proof} 
1. Assume that all the subspaces $A_n$, $n\in\omega$, are winning.
 For every $n\in\omega$ fix a winning strategy
$\Theta_n:\lambda^{<\omega}\to\mathcal P(X)$ of the second player in the
game $CB$ on  $A_n$. Define a strategy
$\Theta:\lambda^{<\omega}\to\mathcal P(X)$ of the second player in the game $CB$
 on  $A=\bigcup_{n\in\omega}A_n$ letting
$\Theta(u_0,\dots,u_n)=\bigcup_{k\le n}\Theta_k(u_k,\dots,u_n)$ for
$(u_0,\dots,u_n)\in\lambda^{<\omega}$. The $u_n$-boundedness of the sets
$\Theta_k(u_k,\dots,u_n)$, $k\le n$, implies the $u_n$-boundedness of
their union $\Theta(u_0,\dots,u_n)$.

We claim that $A\subset\bigcup_{n\in\omega}\Theta(u_0,\dots,u_n)$ for any
infinite sequence $(u_n)_{n\in\omega}\in \lambda^\omega$. Fix any $k\in\omega$.
Regarding the sequence $(u_n)_{n\ge k}$ as the moves of the first
player in the Menger game  on  $A_k$, we see that
$A_k\subset\bigcup_{n\ge k}\Theta_k(u_k,\dots,u_n)$ (according to the
choice of $\Theta_k$ as a winning strategy). Then $$
A=\bigcup_{k\in\omega}A_k\subset \bigcup_{k\in\omega}\bigcup_{n\ge
k}\Theta_k(u_k,\dots,u_n)=\bigcup_{n\in\omega}\bigcup_{k\le
n}\Theta_k(u_k,\dots,u_n)=\bigcup_{n\in\omega}\Theta(u_0,\dots,u_n) $$ and
hence $\Theta$ is a winning strategy of the second player in the  Menger game
 on $A=\bigcup_{n\in\omega}A_n$.

2. Next, assume that all the subspaces $A_n$, $n\in\omega$, are Menger (Hurewicz).
 To show that the union $A=\bigcup_{n\in\omega}A_n$ is
Menger (Hurewicz), fix an infinite sequence of covers
$(u_n)_{n\in\omega}\in\lambda^\omega$. By the Menger (Hurewicz)  property of 
$A_n$, $n\in\omega$, for every $k\in\omega$ there is a ($\gamma$-)cover
$\{B^k_n:n\ge k\}$ of $A_k$ such that each set $B_n^k$, $n\ge k$,
is $u_n$-bounded. Letting $B_n=\bigcup_{k\le n}B_n^k$, we see that
each set $B_n$, $n\in\omega$, is $u_n$-bounded and $\{B_n:n\in\w\}$ is a ($\gamma-$)cover
of $A$.
This proves that the union
$A=\bigcup_{n\in\omega}A_n$ is Menger (Hurewicz).
\end{proof}
Concerning the Scheepers property, the situation with unions is much more delicate.
As it is shown in \cite{BZ}, the class of Scheepers multicovered spaces is closed under
finite unions if and only if two arbitrary ultrafilters are coherent, i.e.
the NCF principle holds, see \cite{Va} for corresponding definitions.

\begin{lemma} \label{l5} 
Let $X$ be a topological space such that the multicovered space $(X,\lambda_{\U(X)})$
is winning (resp. Hurewicz, Scheepers, Menger). Then so is the product
$(X\times Y, \lambda_{\U(X\times Y)})$ for every $\sigma$-compact space $Y$.
\end{lemma}
\begin{proof}
Given arbitrary $\sigma$-compact space $Y$, write it as a union
$\cup\{K_n:n\in\w\}$ of a countable family of its compact subspaces. 
Without loss of generality, $K_n\subset K_{n+1}$ for all $n\in\w$.
Let us denote by $h_n$ the restriction to $X\times K_n$
of the projection $\mathrm{pr}_X:X\times Y\to X$. We claim that $h_n$
is perfect with respect to multicovers  $\lambda_{\U(X\times Y)}|(X\times K_n)$ and $\lambda_{\U(X)}$ 
respectively. Indeed, let $u\in\lambda_{\U(X\times Y)}$ and $d$ be a pseudometric
on $X\times Y$ such that $w=\{B_d(z,1):z\in X\times Y\}$ is inscribed into $u$.
For every $n\in\w$ define a function $d_n:X^2\to\IR$ letting $d_n(x_1,x_2)=\mathrm{sup}\{d((x_1,y), (x_2,y)):y\in K_n\}$
and observe that $d_n$ is a continuous pseudometric on $X$.
Let us fix arbitrary $x\in X$. The perfectness of $h_n$ follows
from $w$-boundedness of $h_n^{-1}(B_{d_n}(x,1/3))$, which can be proven by a
standard argument involving compactness of $K_n$ and the definition of $d_n$. 

Applying Proposition~\ref{basics}(7), we conclude that $(X\times K_n, \lambda_{\U(X\times Y)}|(X\times K_n))$
is winning (resp. Hurewicz, Scheepers, Menger) for all $n\in\w$.
Thus Lemma~\ref{l6} completes our proof in  winning, Hurewicz, and Menger
cases. For the Scheepers property we need some auxiliary arguments.
Assuming that $(X,\lambda_{\U(X)})$ is Scheepers, fix a sequence $(u_n)_{n\in\w}\in\lambda_{\U(X\times Y)}^\w$.
For every $n\in\w$ find $v_n\in\lambda_{\U(X)}$ such that $h_n^{-1}(B)$ is $u_n$-bounded
for every $v_n$-bounded subset $B$ of $X$. Then Lemma~\ref{l3} yields
a proper $\w$-cover $\{B_n:n\in\w\}$ of $X$ such that each $B_n$ is $v_n$-bounded.
It sufficies to show that $\{h_n^{-1}(B_n):n\in\w\}$ is an $\w$-cover of $X\times Y$.
For this purpose fix a finite subset $C=\{(x_i,y_i):i\leq m\}$ of $X\times Y$
and find $n\in\w$ such that $\{x_i:i\leq m\}\subset B_n$ and $\{y_i:i\leq m\}\subset K_n$. 
Then $C\subset B_n\times K_n=h_n^{-1}(B_n)$, which means that $\{h_n^{-1}(B_n):n\in\w\}$
is an $\w$-cover of $X\times Y$ and thus finishes our proof.
\end{proof} 

\noindent \textbf{Remark 1.} 
It is well-known that under additional set-theoretic assumptions
there exists a Hurewicz subspace $S$
of $\IR$ such that $S^2$ is not Menger, see \cite[Theorem 43]{ScTs}.
But this does not  contradict
Lemmas \ref{l2} and \ref{lh}. In order to explain this, let us consider 
Tychonoff spaces $X$ and $Y$. Then the topological space
$X\times Y$ is Menger if and only if so is the multicovered space $(X\times Y,\mathcal O(X\times Y))$,
while the product $(X,\mathcal O(X))\times (Y,\mathcal O(Y))$ is Menger
if and only if so is the multicovered space $(X\times Y,\Op(X)\times\Op(Y))$.
It is easy to see, that $\Op(X)\times\Op(Y)\subset\Op(X\times Y)$,
and these multicovers coincide if  and only if $|X|=1$ or $|Y|=1$,
and consequently $(X\times Y,\mathcal O(X\times Y))$
and $(X\times Y,\mathcal O(X)\times \Op(Y))$ are different multicovered spaces.
But in light of Proposition~\ref{basics}(5) it is more interesting to find out
when the  multicovers $\Op(X)\times\Op(Y)$ and $\Op(X\times Y)$ are isomorphic.
A direct verification shows that this is so when both of them are locally-compact
or Lindel\"of $P$-spaces, while from the above mentioned Theorem~43
of \cite{ScTs}, Proposition~\ref{basics}(5), and Lemma~\ref{lh} we conclude that $\Op(S)^2$ is not equivalent to $\Op(S^2)$.  
\hfill $\Box$
\medskip

\begin{lemma} \label{hur-win-l} 
Let $(X,\lambda)$ be a winning multicovered space. Then there exists a winning strategy 
$\Te_1$ of the second player in the game $CB$ such that $\{\Te_1(u_0,\ldots,u_{n}):n\in\w\}$
is a $\gamma$-cover of $X$ for all sequences $(u_n)_{n\in\w}\in\lambda^\w$.
 \end{lemma}
\begin{proof}
Let us fix some
  winning strategy $\Theta$
of the second player in the game $CB$ on $(X,\lambda)$.

We claim that the map $\Theta_1:\lambda^{<\omega}\rightarrow\mathcal{P}(X)$,
\[ \Theta_1:(u_0,u_1,\ldots,u_n)\mapsto\bigcup_{ 0\leq i_0\leq i_1\leq\cdots\leq i_k=n } \Theta(u_{i_0},u_{i_1},\ldots,u_{i_k}) \]
is a winning strategy of the second player in the  game $CB$ on $(X,\lambda)$
with the required property.

Suppose, to the contrary, that there exists a sequence of covers $(u_n)_{n\in\omega}\in\lambda^{\omega}$,
a subsequence $(i_k)_{k\in\omega}\in\omega^{\omega}$, and $x\in X$  
such that
$x\not\in\bigcup_{k\in\omega}\Theta_1(u_0,u_1,\ldots,u_{i_k})$. But $\Theta$
is a winning strategy in the Menger game on $(X,\lambda)$, which together with definition of $\Theta_1$
gives us 
$X=\bigcup_{k\in\omega}\Theta(u_{i_0},u_{i_1},\ldots,u_{i_k})\subset\bigcup_{k\in\omega}\Theta_1(u_0,u_1,\ldots,u_{i_k})$,
a contradiction.
\end{proof}

\begin{corollary}\label{cb1} 
The class of winning multicovered spaces is closed under finite products of its
elements.
\end{corollary}
\begin{proof}
Let $(X,\lambda)$ and $(Y,\nu)$ be two winning multicovered spaces and $\Te_X$ and
$\Te_Y$ be winning strategies of the second player in the game $CB$ on $(X,\lambda)$
and $(Y,\nu)$ respectively having the property from Lemma~\ref{hur-win-l}.
A direct verification shows  that the strategy 
\[  \Te:(u_0\cdot v_0,\ldots,u_n\cdot v_n)\mapsto\Te_X(u_0,\ldots,u_n)\times\Te_Y(v_0,\ldots,v_n)   \] 
is winning in the game $CB$ on the product $(X\times Y,\lambda_X\times\lambda_Y)$.
\end{proof}

The next corollary answers \cite[Problem 1]{HRT} in negative.

\begin{corollary}\label{Banakh} 
The product of finitely many strictly $o$-bounded topological groups is strictly
$o$-bounded.
\end{corollary}
\begin{proof}
Follows from the observation that a topological group $G$ is strictly $o$-bounded
if and only if the multicovered space $(G,\lambda_L(G))$ is winning, see Corollary~\ref{cb1},
and Proposition~\ref{basics}(2,5).
\end{proof}

The following lemma is the central part of the work.
\begin{lemma} \label{newcrucial} 
 Let $G$ be a topological group and $X\subset G$ be a set of its generators.
If the multicovered space $(X\cup X^{-1},\lambda_R(G)|X\cup X^{-1})$ is winning, then so is
$(G,\lambda_R(G))$. 

If, additionaly, $G$ is abelian and $(X,\lambda(G)|X)$
is Hurewicz (Scheepers), then so is $(G,\lambda(G))$.
\end{lemma}
\begin{proof}
1. We start by proving the ``winning'' part.
Assuming that $(X\cup X^{-1},\lambda_R(G)|X\cup X^{-1})$ is winning, find a strategy 
$\Te:\lambda_R^{<\w}\rightarrow\mathcal P(G)$ such that $\{\Te(u_0,\ldots, u_k):k\in\w\}$
is a $\gamma$-cover of $X\cup X^{-1}$ for every sequence $(u_n)_{n\in\w}\in\lambda_R^\w$.
Let $\mathcal B$ be the family of all open neighborhoods of the identity of $G$.

Next, for every $s\in\lambda_R^{<\w}$ we shall construct a sequence $w(s)=(w(s)_n)_{n\in\w}\in\lambda_R^\w$.  
Let $s=(u_0,\ldots,u_m)$,  $U\in\mathcal B$ be such that $u_m\prec\{Uz:z\in G\}$,
and $U_0\in\mathcal B$ be such that $U\supset U_0^2$.
Put $w(s)_0=\{U_0 z:z\in G\}$ and $A_0(s)=\Te(s\hat{\ }w(s)_0)$. Assume that for some $n\in\w$
and for all $k\leq n$ we have already constructed $w(s)_k=\{U_k z:z\in G\}\in\lambda_R$ and
$A_k(s)\subset G$ such that the following conditions are satisfied:
\begin{itemize}
\item[$(i)$] $A_k(s)=\Te(s\hat{\ }w(s)_0\hat{\ }\cdots\hat{\ }w(s)_k)$;
\item[$(ii)$] $U_k\supset U_l^2$ for all $k<l\leq n$;
\item[$(iii)$] $A_k(s) B$ is $w(s)_{k-1}$-bounded for every $w(s)_l$-bounded subset
$B$ of $G$, where $k<l\leq n$ and $w(s)_{-1}=u_m$.
\end{itemize} 
Since $A_n(s)$ is $\{U_n z:z\in G\}$-bounded, there exists a finite subset
$K$ of $G$ such that $A_n(s)\subset U_nK$. Let us find
$U_{n+1}\in\mathcal B$  such that $zU_{n+1}z^{-1}\subset U_{n}$ for all $z\in K$
and $U_{n+1}^2\subset U_n$, and set $w(s)_{n+1}=\{U_{n+1}z:z\in G\}$.
Given any $k<n+1$ and an $w(s)_{n+1}$-bounded subset $B$ of $G$, consider the product
$C=A_k(s)B$. If $k<n$, then the $w(s)_{k-1}$-boundedness of $C$ follows from 
$(iii)$ and the equation $w(s)_n\prec w(s)_{n+1}$. Thus, it sufficies to consider the case
$k=n$. Let $L$ be a finite subset of $G$ such that $B\subset U_{n+1}L$.
Then 
$$C= A_n(s)B \subset U_n K U_{n+1}L\subset U_nU_{n}KL\subset U_n^2KL\subset U_{n-1}KL,$$   
which yields the $w(s)_{n-1}$-boundedness of $C$, and thus completes our inductive construction
of the sequence $(w(s)_n)_{n\in\w}$ satisfying $(i)-(iii)$ for all $n\in\w$.
Observe, that the condition $(iii)$ implies that the product $A_0(s)A_2(s)\cdots A_{2n}(s)$
is $u_m=w(s)_{-1}$-bounded for all $n\in\w$. 

Given any $s=(u_0,\ldots, u_{n-1}) \in\lambda_R^{<\w}$,
construct a finite sequence $(q_0(s),\ldots,q_{2n-2}(s)) \in (\lambda_R^{<\w})^{<\w}$
as follows:
$$q_0(s)=(u_0),\  q_{2k+1}(s)=q_{2k}(s)\hat{\ }w(q_{2k}(s))|(2k+1),\ q_{2k+2}(s)=q_{2k+1}(s)\hat{\ }u_{k+1}.$$
Let $\Te_1(s)=A_0(q_{2n-2}(s))A_2(q_{2n-2}(s))\cdots A_{2n-2}(q_{2n-2}(s))$. We claim that 
$\Te_1$ is a winning strategy of the second player in the game $CB$  on
$(G,\lambda_R)$. Indeed, from the above it follows that $\Te_1(s)$ is 
$w_{-1}(q_{2n-2}(s))=u_{n-1}$-bounded,
which implies that $\Te_1$ is a strategy of the second player.
To show that it is winning, consider arbitrary $z=x_0 x_1\cdots x_m\in G$, where
$x_i\in X\cup X^{-1}$ for all $i\leq m$. Let $t=(u_n)_{n\in\w}\in\lambda_R^\w $ be a sequence of
covers of $G$. Our proof will be completed as soon as we show 
that there exists $n\in\w$ such that $\Te_1(t|n)\ni z$.
For this aim consider the sequence $(v_k)_{k\in\w}\in\lambda_R^\w$ such that
for every $n\in\w$ there exists $k_n\in\w$ such that $q_{2n-2}(t|n)=(v_{0},\ldots,v_{k_{n-1}})$
(the definition of $q_{-}(-)$ easily yields such a sequence, and $k_n=k_{n-1}+(2n-2)+1$).
From the above it follows that 
\begin{eqnarray*} \Te_1(t|n)=A_0(q_{2n-2}(t|n)) A_2(q_{2n-2}(t|n))\cdots A_{2n-2}(q_{2n-2}(t|n))=\ \ \ \ \ \ \ \ \ \ \ \ \ \ \ \ \ \ \ \ \ \ \ \   \\
= \Te(v_0,\ldots,v_{k_{n-1}}) \Te(v_0,\ldots,v_{k_n}, v_{k_{n-1}+1}, v_{k_{n-1}+2})\cdots \Te(v_0,\ldots,v_{k_{n-1}},\ldots, v_{k_{n-1}+2n-2}). 
\end{eqnarray*}
By our choice of $\Te$, the family $\{\Te(v_0,\ldots,v_k):k\in\w\}$ is a $\gamma$-cover of $X\cup X^{-1}$,
consequently there exists $l\in\w$ such that $\{x_0,\ldots,x_m\}\subset\Te(v_0,\ldots,v_k)$
for all $k\geq l$. Let $n>m$ be such that $k_{n-1}>l$.
Then $\{x_0,\ldots, x_m\}\subset\Te(v_0,\ldots,v_{k_{n-1}},\ldots,v_{k_{n-1}+2i})$ for all
$i\in\{0,2,\ldots,2n-2\}$, which implies $z\in\Te_1(t|n)$.

2. Let us assume that the multicovered space $(X,\lambda(G)|X)$
is Scheepers and set $\lambda=\lambda(G)$.
Given a sequence
 $(u_n)_{n\in\w}\in\lambda^\w$, find a sequence $(O_n)_{n\in\w}$
of open neighborhoods of the neutral element $e$ such that
$u_n\prec\{g+O_n : g\in G\}$, $-O_n=O_n$, and $2 O_{n+1}\subset O_n$ for all $n\in\w$.
By the definition of the Scheepers property applied to $(X,\lambda|X)$ there
exists a sequence $(K_n)_{n\in\w}$ of finite subsets of $G$ such that
the family
$v=\{K_n+O_n : n\in\w\}$ is a proper $\w$-cover of $X$.
Without loss of generality, $K_n=-K_n$ and $K_n+K_n\subset K_{n+1}$ for all
$n\in\w$. 
We claim that $v_1=\{K_{2n}+O_n : n\in\w\}$
is an $\w$-cover of $G$.

Indeed, from the above it
follows that $K+O_n\supset X$ for every $n\in\w$, where $K=\bigcup_{n\in\w}K_n$.
Consequently for every $x\in X$ we can define a nondecreasing number sequence $z(x)$
letting $z(x)_n=\min\{m\in\w:x\in K_m+O_n\}$.
Since $v$ is a proper $\w$-cover of $X$, for every finite subset $S$ of $X$
the set $I_S=\{n\in\w:z(x)_n\leq n \mbox{ for all }x\in S\}$ is infinite.
Now, consider  arbitrary finite subset $A$ of $G$ and find some finite 
subset $S$ of $X$ and $m\in\IN$ such that $A\subset m(S-S)$. Let us fix arbitrary
$l\in I_S\cap [3m,+\infty)$. Then $S-S\subset K_l+O_l-K_l-O_l\subset K_{l+1}+O_{l-1}$,
$2(S-S)\subset 2(K_{l+1}+O_{l-1})\subset K_{l+2}+O_{l-2}$ and so on.
Proceeding in this fashion, we obtain $A\subset m(S-S)\subset K_{l+m}+O_{l-m}$. 
Since $l\geq 3m$, there exists $n\in\w$ such that  $n\leq l-m$
and $2n\geq l+m$, which yields $O_n\supset O_{l-m}$ and $K_{2n}\supset K_{l+m}$.
From the above it follows that $K_{l+m}+O_{l-m}\subset K_{2n}+O_n$,
which proves that the multicovered space $(G,\lambda)$ is Scheepers.

3. The proof of the ``Hurewicz'' part is similar to that of the ``Scheepers'' one.
Let   $\lambda$,
 $(u_n)_{n\in\w}\in\lambda^\w$, and $(O_n)_{n\in\w}$
be such as in the previous item.
Since the multicovered space $(X,\lambda|X)$ is Hurewicz,  there
exists a sequence $(K_n)_{n\in\w}$ of finite subsets of $G$ such that
the family
$v=\{K_n+O_n : n\in\w\}$ is a $\gamma$-cover of $X$.
Without loss of generality, $K_n=-K_n$ and $K_n+K_n\subset K_{n+1}$ for all
$n\in\w$. 
We claim that $v_1=\{K_{2n}+O_n : n\in\w\}$
is a $\gamma$-cover of $G$.
Let $K$ and $z(x)\in\w^{\uparrow\w}$ be such as in the second item.
Since $v$ is a  $\gamma$-cover of $X$, for every finite subset $S$ of $X$
the set $I_S=\{n\in\w:z(x)_n\leq n \mbox{ for all }x\in S\}$ is cofinite, i.e.
the complement $\w\setminus I_S$ is finite.
Now, consider  arbitrary $z\in G$ and find some finite 
subset $S$ of $X$ and $m\in\IN$ such that $z\in m(S-S)$. Let us fix some
$l\geq 3m$ such that $[l,+\infty)\subset I_S$. Then for every $p\geq l$ we have $S-S\subset K_p+O_p-K_p-O_p\subset K_{p+1}+O_{p-1}$,
$2(S-S)\subset 2(K_{p+1}+O_{p-1})\subset K_{p+2}+O_{p-2}$ and so on.
Proceeding in this fashion, we obtain $z\in m(S-S)\subset K_{p+m}+O_{p-m}$. 
Since $p\geq l\geq 3m$,  $p+m\leq 2(p-m)$, and consequently   $K_{p+m}\subset K_{2(p-m)}$,
which yields $z\in K_{2(p-m)}+O_{p-m}$.
Since $p\geq l$ was chosen arbitrary, we conclude that $z\in O_n+K_{2n}$ for all
$n\geq l-m$, which means that $v_1$ is a $\gamma$-cover of $G$.
\end{proof}
\noindent\textbf{Remark 2.}
The winning property of any  abelian topological group $H$ containing $X$
as a set of its generators can be derived from the winning property of $(X,\lambda(H)|X)$ 
much easier than in the  general  case considered in Lemma~\ref{newcrucial}.
Given  any finite sequence $j=(j_0,\ldots,j_{n-1})\in \{-1,1\}^{<\w}$,
define a map $\psi_j:X^n\to G$  letting $\psi_j(x_0,\ldots,x_{n-1})=j_0x_0+\cdots+j_{n-1}x_{n-1}$.
A direct verification shows that $\psi_j$ is uniformly-bounded with respect to multicovers
$(\lambda(H)|X)^n$ and $\lambda(H)$ (here commutativity is essentially used), and hence 
$(\psi_j(X^n), \lambda(H)|\psi_j(X^n))$ is winning for each $j\in\{-1,1\}^{<\w}$.
Now it sufficies to use Corollary~\ref{cb1}, Proposition~\ref{basics}(8), and Lemma~\ref{l6}.

The same arguments work for the Hurewicz property. In case of the Scheepers property
we have to additionaly prove that the countable union of uniformly-bounded
images of finite powers of a Scheepers space $(X,\lambda(H)|X)$  is Scheepers
(the union of Scheepers multicovered spaces could be not Scheepers, see the discussion following Lemma~\ref{l6}).
\hfill $\Box$
\medskip

\noindent \textbf{Proofs of Theorems~\ref{main1} and \ref{main-mult1}.}
We shall prove these theorems by showing that the conditions $(1)-(8)$
of Theorem~\ref{main1} are equivalent to the Scheepers property of $(X,\lambda_{\U(X)})$
(note that Theorem~\ref{main-mult1} states that $(1)$ is equivalent to the Scheepers property
of $(X,\lambda_{\U(X)})$),
and the last condition will be denoted by $(9)$.
The implication $(2)\Rightarrow(1)$ is obvious.
The implications $(5)\Rightarrow (4)$ and $(4)\Rightarrow(3)$
follow from the continuity of linear maps $\varphi:\ls(X)\to L(X)$
and $\psi:L(X)\to L_p(X)$ extending the identity map $\mathrm{id}_X$,
and the simple fact  that the $o$-boundedness
is preserved by continuous homomorphic images, see, e.g., \cite{Tk98}. 

In addition, we shall prove the subsequent implications: 
$(1)\Rightarrow(9)$, $(9)\Rightarrow(2)$, $(9)\Rightarrow(5)$, $(3)\Rightarrow(9)$,
$(6)\Leftrightarrow (9)$, $(7)\Leftrightarrow (1)$,  $(8)\Leftrightarrow (1)$.

$(1)\Rightarrow(9)$.  Since $A(X)$ is $o$-bounded, it is $\w$-bounded,
and thus the uniform space $(X,\U(X))$ as well as the multicovered space
$(X,\lambda_{\U(X)})$ are $\w$-bounded by Proposition~\ref{fromG}. Therefore
$X$ and $G=A(X)$ satisfy the conditions of Lemma~\ref{l1}, and consequently for every
$n\in\IN$ the map $\psi_n$ defined there is perfect with respect
to $\lambda_{\U(X)}^n$ and $\lambda(A(X))$. As it was stressed in Introduction,
the $o$-boundedness of the group $A(X)$ is equivalent to the Menger property
of the multicovered space $(A(X),\lambda(A(X)))$. Applying Proposition~\ref{basics}(7),
we conclude that $(X^n,\lambda_{\U(X)}^n)$ is Menger for all $n\in\IN$,
and consequently $(X,\lambda_{\U(X)})$ is Scheepers by Lemma~\ref{l2}.

$(9)\Rightarrow(2)$. Assume that $(X,\lambda_{\U(X)})$ is Scheepers. 
Then so is $(X,\lambda(A(X))|X)$. Applying Lemma~\ref{newcrucial}, we conclude
that $(A(X),\lambda(A(X)))$ is Scheepers too, and thus $(A(X)^n,\lambda(A(X))^n)$
is Menger for all $n\in\IN$ by Lemma~\ref{l2}, which means that $A(X)^n$
is $o$-bounded for all $n\in\IN$.

$(9)\Rightarrow(5)$. Let us note, that we have already proven the equivalence
of items $(1)$, $(2)$, and $(9)$. Let $X$ be a topological space satisfying $(9)$.
Then $(X,\lambda_{\U(X)})$ is Scheepers, and hence so is the multicovered space
$(X\times\IR,\lambda_{\U(X\times\IR)})$ by Lemma~\ref{l5}. Consequently
$A(X\times \IR)^n$ is $o$-bounded for all $n\in\IN$.
Consider a map $h:X\times \IR\to \ls(X)$, $h(x,r)=rx$. Since $\ls(X)$
is a linear topological space, $h$ is continuous, and hence it admits a continuous extension
to a homomorphism $\tilde{h}:A(X\times\IR)\to \ls(X)$. A direct verification shows that
$\tilde{h}$ is surjective. Therefore  
$\ls(X)$ is a continuous homomrphic image of $A(X\times\IR)$,
and consequently $\ls(X)^n$  is a continuous homomorphic image
of $A(X\times\IR)^n$ for all $n\in\IN$. 
 
$(3)\Rightarrow(9)$. It sufficies to use the fact that $X$ and $G=L_p(X)$
satisfy the conditions of Lemma~\ref{l1}, see \cite[Chapter~0]{Ar89},
and apply the same argument as in the proof of the implication
$(1)\Rightarrow (9)$.

$(9)\Leftrightarrow(6)$. Assuming that $(X,\lambda_{\U(X)})$  is Scheepers,
fix a continuous surjective function $f:X\to Y$ onto a metrizable space $Y$.
Then $f$ is uniformly bounded with respect to multicovers
$\lambda_{\U(X)}$ and $\lambda_{\U(Y)}$ being uniformly continuous
with respect to uniformities $\U(X)$ and $\U(Y)$.
Therefore $(Y,\lambda_{\U(Y)})$ is Scheepers by Proposition~\ref{basics}(8).
In particular, this implies that $Y$ is Lindel\"off and hence $\lambda_{\U(Y)}$
 is equivalent to $\Op(Y)$ by Corollary~\ref{cle213}. Then the multicovered space
$(Y,\Op(Y))$  is Scheepers by Proposition~\ref{basics}(5).

Next, assume that $(X,\lambda_{\U(X)})$ is not Scheepers and  find
a sequence $(u_n)_{n\in\w}\in\lambda_{\U(X)}^\w$ witnessing for this.
For every $n$ find an entourage $U_n\in\U(X)$ such that $w_n=\{U_n(x):x\in X\}$
is inscribed into $u_n$. Let $d$ be a continuous pseudometric on $X$
such that $v_n=\{B_d(x,2^{-n}):x\in X\}$ is inscribed into $w_n$.
Then the identity map $\mathrm{id}_X$ is perfect with respect to multicovers
$\lambda_1=\{u_n:n\in\w\}$ and $\lambda_2=\{v_n:n\in\w\}$ of $X$.
Since $(X,\lambda_1)$ is not Scheepers, so is $(X,\lambda_2)$
 by Proposition~\ref{basics}(7). Consequently $X$ endowed with the topology
generated by $d$ is not Scheepers, and hence there are non-Scheepers metrizable images
of $X$.
 
$(7)\Leftrightarrow(1)$. Let us note that we have already proven the equivalence
$(1)\Leftrightarrow(6)$.
It is well known that every Lindel\"of space is Hewitt-complete and  every
continuous map $f:X\to Y$ from a space $X$ into a Hewitt-complete space
$Y$ extends to a continuous map $\hat{f}:\nu X\to Y$, see \cite[Th.~3.11.12, 3.11.16]{En}. 
Let $X$ be such that $A(X)$ is $o$-bounded and $Y$ be a continuous metrizable image
of $\nu X$ under a map $f$.
Then $Y$ is Lindel\"of containing a dense  Lindel\"of (even Scheepers) subspace  $Z=f(X)$. 
Therefore $Y$ as well as $Z$  are Hewitt-complete, and hence the map
$f|X$ extends  to a continuous map $g:\nu X\to Z$. Since $f$ and $g$ coincide
on the dense subset $X$ of $\nu X$, we get $f=g$, and hence $Y=Z=f(X)$.
Thus we have already proven that each continuous metrizable image of $\nu X$
is Scheepers, which implies the $o$-boundedness of $A(\nu X)$.

Now, assume that $A(\nu X)$ is $o$-bounded. It follows that each metrizable image
of $\nu X$ is Scheepers. The same argument as in the previous paragraph gives 
that each metrizable image of $X$ is Scheepers as well, and hence $A(X)$
is $o$-bounded.

$(8)\Leftrightarrow (1)$.
It is well-known  \cite[8.5.8(b)]{En}
that there are natural embeddings of $\nu X$ and $\mu X$ into the Stone-\v{C}ech
compactification $\beta X$ such that $X\subset\mu X\subset\nu X\subset\beta X$.
This permits us to apply the same argument as in the proof of the equivalence
of $(1)$ and $(7)$ and conclude that  $X$ and $\mu X$ have the same continuous metrizable images,
and then apply  already proven equivalence $(1)\Leftrightarrow(6)$.
\hfill $\Box$

\medskip

\noindent \textbf{Proof of Theorem~\ref{main-win}.}
A part of the  proof of this theorem runs fairly in a similar way 
as that of Theorem~\ref{main1}.  Namely,
the implications 
$(4)\Rightarrow(8)$, $(8)\Rightarrow(2)$,
$(4)\Rightarrow(7)$, and $(5)\Rightarrow(8)$
 can be proven similarly to the implications
$(1)\Rightarrow(9)$, $(9)\Rightarrow(2)$, 
$(2)\Rightarrow(5)$, $(3)\Rightarrow(9)$
of Theorem~\ref{main1} respectively (one has to additionaly use that
the product of winning multicovered spaces is winning, and $(X\cup X^{-1},\lambda_R(X\cup X^{-1}))$
is winning provided so is $(X,\lambda_{\U(X)})$ by Lemma~\ref{l6}).

The implications $(3)\Rightarrow(2)$, $(2)\Rightarrow(1)$, and $(7)\Rightarrow(6)$,
$(6)\Rightarrow(5)$
immediately follow from corresponding definitions. Thus we are left with the
task of proving the implications $(1)\Rightarrow(4)$ and $(2)\Rightarrow(3)$.
Concerning the implication $(1)\Rightarrow(4)$, it follows from Proposition~\ref{basics}(8)
and the fact that
the homomorphism $f:F(X)\to A(X)$ extending the identity map on $X$
is uniformly continuous with respect to uniformities $\U_{L\wedge R}(F(X))$
and $\U(A(X))$, and hence is uniformly bounded with respect to multicovers
$\lambda_{L\wedge R}(F(X))$ and $\lambda(A(X))$ of $F(X)$ and $A(X)$
respectively. 

$(2)\Rightarrow(3)$. Let us note, that in light of Corollary~\ref{cb1}
it sufficies to prove that the multicovered space $(F(X),\lambda_{L\vee R})$ 
is winning. Set  $\Delta_{F(X)}=\{(x,x):x\in F(X)\}$. 
Then the map  $i:X\ni x\mapsto (x,x)\in\Delta_{F(X)}$ as well as its inverse
are obviously perfect with respect to multicovers
 $\lambda_{L\vee R}$ and  $\lambda_L\times \lambda_R|\Delta_{F(X)}$ of $F(X)$ and $\Delta_{F(X)}$ respectively.
Since $F(X)$ is strictly $o$-bounded, both of the multicovered space $(F(X),\lambda_R)$
and $(F(X),\lambda_L)$ are winning, and hence so is the product
$(F(X)^2,\lambda_R\times\lambda_L)$, and finally the multicovered spaces
$(\Delta_{F(X)},\lambda_R\times\lambda_L|\Delta_{F(X)})$ and  $(F(X),\lambda_{L\vee R})$
are winning as well.
\hfill $\Box$
\medskip

\begin{example} \label{count-com}
There exists a countably-compact space $Z$ and a pseudometric $d$
on $Z^2$ such that the pseudometric space $(Z^2,d)$ is not Lindel\"of.
 \end{example}
\begin{proof}
To begin with, let us note that it sufficies to construct two countably-compact
spaces $X$ and $Y$ and a pseudometric $d$ on their product $X\times Y$
such that the corresponding pseudometric space is not Lindel\"of,
and then the topological sum $Z$  of $X$ and $Y$ obviously admits a required
pseudometric.
Let $D$ be a discrete space of size $|D|=\aleph_1$.
 Similarly to Example~3.10.19 of \cite{En} 
 we define a function $f$ assigning to each countable subset  $A$ of $\beta D$  
some element $f(A)\in\overline{A}\setminus A$. Let $X_0=D$
and $$ X_\alpha=(\bigcup_{\gamma<\alpha}X_\gamma)\cup f([\bigcup_{\gamma<\alpha}X_\gamma]^{\aleph_0})$$
for $0<\alpha<\w_1$, where $[A]^{\aleph_0}$ stands for the family of all countable subsets
of a set $A$. Thus we have already defined  a transfinite sequence $(X_\alpha)_{\alpha<\w_1}$
of subsets of $\beta D$. The space $X=\bigcup_{\alpha<\w_1}X_\alpha$ is obviously 
countably-compact
(every countable subset $A$ of $X$ is contained in some $X_\alpha$, and hence
is not closed in $X$). It is easy to prove that $|X|\leq\mathfrak c$. Set
$Y= D\cup(\beta D\setminus X)$. According to Theorem~3.6.14 of \cite{En},
 $|\overline{A}|=2^{\mathfrak c}$ for every
countable $A\subset\beta D$, and hence $Y$ is countably-compact as well.
It sufficies to observe that $X\times Y$ contains an open discrete subspace
$\Delta_{D}=\{(x,x):x\in D\}$ of size $\aleph_1$, and hence admits a
non-Lindel\"of pseudometric.  
\end{proof}

\noindent \textbf{Proof of Theorem~\ref{main-hur}.}
The proof of this this theorem is quite similar to that of Theorem~\ref{main1}
and is left to the reader.
\hfill $\Box$
\medskip
 
\noindent \textbf{Remark 3.}
The characterization of spaces $X$ such that $F(X)$ is $o$-bounded   is the same as
in the abelian case. But its proof requires a technique of (semi)filter games
investigated by C.~Laflamme,
and is not within the methods used here. 
This problem is to be considered in \cite{BZ} from a more general point of view. 
\hfill $\Box$
\smallskip

\noindent \textbf{Acknowledgments.} The author wishes to express his thanks
to  Taras Banakh and  Igor Guran for many stimulating conversations.

\noindent Department of Mechanics and Mathematics, \\
 Ivan Franko Lviv National University, \\
Universytetska 1, Lviv, 79000, Ukraine.

\textit{E-mail address:}   \texttt{lzdomsky@rambler.ru}

\end{document}